\documentclass[12pt]{article}

\usepackage{sola}

\function{Gma}{\Gamma}

\defoperator{ii}{I}

\defoperator{QA}{Q^A}
\defoperator{QL}{Q^L}
\defoperator{QB}{Q^B}
\defoperator{QF}{Q^F}
\def\L{{\cal L}}
\def\P{{\cal P}}

\defoperator{MA}{{\cal M}^\star}

\def\iig_#1[#2]{I_{#1}\br[{#2}]}
\defoperator{iigg}{I_g}
\def\iigf[#1]{\iigg[f,#1]}

\def\QAg_#1[#2]{Q^A_{#1}\br[{#2}]}
\defoperator{QAgg}{Q^A_g}
\def\QAgf[#1]{\QAgg[f,#1]}

\def\QLg_#1[#2]{Q^L_{#1}\br[{#2}]}
\defoperator{QLgg}{Q^L_g}
\def\QLgf[#1]{\QLgg[f,#1]}

\def\QBg_#1[#2]{Q^B_{#1}\br[{#2}]}
\defoperator{QBgg}{Q^B_g}
\def\QBgf[#1]{\QBgg[f,#1]}

\def\QFg_#1[#2]{Q^F_{#1}\br[{#2}]}
\defoperator{QFgg}{Q^F_g}
\def\QFgf[#1]{\QFgg[f,#1]}

\function{Lv}{\L[v]}

\def\sig_#1(#2){\sigma_{#1}\fpr(#2)}

\def\ewp(#1){\E^{\I \omega \pr({#1})}}
\def\ewg(#1){\E^{\I \omega g\pr(#1)}}

\def\sjacobiandd_#1^#2{\vc{J}_{#1}^{#2}}
\def\sjacobiand_#1^#2{\sjacobiandd_{#1}^{#2}(\xx)}

\def\sjacobian_#1{\sjacobiand_{#1}^d}

\def\fO(#1){\Oh\fpr({#1})}

\def\Leqpn^#1#2{\L[v]^{(#1)}(#2)=f^{(#1)}(#2)}
\def\Leqnpn^#1#2{\Leqpn^{#1}{x_#2}}

\def\xx{\vc x}

\def\Km_#1{{\cal K}_{#1}}
\defoperator{K}{{\cal K}_n}

\begin{document}

\setupproofs
\setupenvironments

\def\H{{\cal H}}
\def\Mobius{M\"obius}
\def\Holder{H\"older}

\def\soitem{\newitem}

\def\sociteb[#1;#2]{\cite[#2]{#1}}
\def\socite#1{\sociteb[#1]}

\refmodetrue

\titled={Change of variable formul\ae\ for regularizing  slowly decaying and oscillatory Cauchy and Hilbert transforms }


\maketitle

\Abstract
	Formul\ae\ are derived for expressing Cauchy and Hilbert transforms of a function $f$ in terms of Cauchy and Hilbert transforms of $f(x^r)$.
	 When $r$ is an integer, this corresponds to evaluating the Cauchy transform of $f(x^r)$ at all choices of $z^{1/r}$.     Related formu\ae\ for rational $r$ result in a reduction to a generalized Cauchy transform living on a Riemann surface, which in turn is reducible to the standard Cauchy transform.  
	 These formul\ae\ are  used to regularize the behaviour of functions that are slowly decaying or oscillatory, in order to facilitate numerical computation and extend asymptotic results.
	
\Section{intro} Introduction.

	The (one-sided) \newterm{Cauchy transform}  is defined for $z \notin [0,\infty)$ by
	$$\CC f(z) =  {1 \over 2 \pi \I} \int_0^\infty {f(t) \over t - z}\dx.$$
A related transform is  the (one-sided) \newterm{Hilbert transform} defined for $x \in (0,\infty)$ by
	$${\cal H} f(x) = {1 \over \pi} \dashint_0^\infty {f(t) \over t - x} \dt = {1 \over \pi} \lim_{\epsilon \rightarrow 0^+} \pr(\int_0^{x-\epsilon} + \int_{x+\epsilon}^\infty) {f(t) \over t - x } \dt,$$
i.e.,  the integral is a \newterm{Cauchy principle-value integral}.   Applications of Cauchy and Hilbert transforms are numerous, including integrable systems \cite{AblowitzInverseScattering,AblowitzClarksonSolitons}, random matrix theory \cite{DeiftOrthogonalPolynomials},  and signal processing \cite{HahnHilbert,KingHilbertTransforms1,KingHilbertTransforms2}.

In this paper, Cauchy and Hilbert transforms of $f$ are expressed in terms of Cauchy and Hilbert transforms of $f(x^r)$.  Cauchy and Hilbert transforms are classical objects that have been thoroughly studied, and  their behaviour under many transformations  are known (cf.~\socite{KingHilbertTransforms2; Appendix 1} and \socite{TheCauchyTransform; Chapter 11}).  Suprisingly,   formul\ae\ for these transforms under the  change of variables $x^r$  appear to be unknown.

The following  three computational applications motivate the results:
\beginorderedlist
	\soitem If $f$ decays, then $f(x^r)$ will have an increased rate of decay for large $r$, and the domain of integration can be truncated.
	\soitem If $f$ has a slowly decaying fractional asymptotic series at  infinity
	$$f(x) \sim \sum_{k = 0}^\infty f_k x^{k \over r},$$
then $f(x^r)$ has a standard asymptotic series, and the numerical method of \cite{SOHilbertTransform} can be employed.   Hilbert transforms of slowly decaying functions have applications to the critical behaviour of nonlinear PDEs \cite{AntunesNonlinearPDE}.      
	\soitem Even for small $r$, $f(x^r)$ may have a special form that makes its Cauchy and Hilbert transforms amenable to approximation.  For example, if $f$ is highly oscillatory with a high-order oscillator
	$$f(x) = \E^{\I \omega x^{1/r}} g(x),$$
then $f(x^r)$ has a standard Fourier oscillator.  Hilbert transforms with the   standard Fourier oscillator have been investigated asymptotically \cite{WongHilbertTransform,WangOscillatoryHilbert} and numerical methods have been developed  \cite{KellerOscillatoryHilbert,WangOscillatoryHilbert}.
 Cauchy transforms of oscillatory functions underly  Riemann--Hilbert problems  throughout integrable systems  \cite{AblowitzInverseScattering,AblowitzClarksonSolitons}, random matrix theory \cite{DeiftOrthogonalPolynomials}, and have applications to  gravity waves  \cite{ChapmanGravityWaves}.  
\endorderedlist

	The paper is structured as follows. We first review Plemelj's lemma  in \secref{caucPlem}, which underlies subsequent proofs.  We   consider integer $r$ in \secref{map}, in which case we can express the Cauchy transform of $f$ in terms of the Cauchy transform of $f(x^r)$ evaluated at {\it all} $r$ choices of $z^{1/r}$, see \lmref{Pol}.    The case where $r = 1/q$ and $q$ is an integer is discussed in \secref{inte},  obtaining a representation in terms of a generalized Cauchy transform of $f(x^r)$ living on a $q$-sheeted Riemann surface, see \lmref{int}.   The two preceding versions are combined for $r = p/q$, in which case we need to evaluate a generalized  Cauchy transform of $f(x^r)$ living on a $q$-sheeted Riemann surface at $p$ different points, see \thref{rational}.  	

	 We  turn our attention to examples in \secref{num}, exploiting the new representations for each of the above motivating cases.  We compare two approaches for numerically calculating slowly decaying Cauchy transforms: mapping the Cauchy transform to $f(x^r)$ for large $r$ to induce fast decay,  and choosing $r$ so that $f(x^r)$  is in a form amenable to calculation via \cite{SOHilbertTransform}.     We also consider asymptotic expansions of oscillatory Hilbert transforms,  combining the proposed formul\ae\ with the asymptotic expansion derived in \cite{WangOscillatoryHilbert}.

 In \appref{irr}, we   propose a candidate  formula relating the Cauchy transform of $f$ to $f(x^r)$ for  irrational $r$, which is given in terms of evaluating a Cauchy transform living on an infinite sheeted Riemann surface at an infinite number of points.  This formula is not proved, but numerical evidence support its validity.  Unfortunately, the formula has poor convergence properties, and hence is of limited practical use. 

\Remark
	The function $z^r$ will always take the standard principle branch with a branch cut along $(-\infty,0]$.

\Section{caucPlem} Cauchy transform, Hilbert transform and Plemelj's lemma.

Linking the Cauchy transform and Hilbert transform, and playing a crucial role in the theory below, is Plemelj's lemma (alternatively known as Sokhotskyi formulas).  We first define a (sufficient) class of functions so that Plemelj's lemma applies:

\Definition{unifHolder}
	A function $f$ is \Holder-continuous at $x$ if there exists positive $\gamma$ and $\mu$ such that
	$$\abs{f(x+h) - f(x)} \leq \mu \abs{h}^\gamma$$
for sufficiently small $h$.  A function $f$ is uniformly \Holder-continuous on $[a,b]$ if it is \Holder-continous for every $x \in [a,b]$ with the same choice of $\gamma$ and $\mu$. 
Finally,  $f$  is \Holder-continuous on $(0,\infty)$ if it is uniformly \Holder-continuous on every bounded subset $[a,b] \subset (0,\infty)$ and is \Holder-continuous at infinity:
	$f(1/x)$ is \Holder-continuous at 0.  

Plemelj's lemma states the following:

\Lemma{plemelj}
	Suppose $f$ is  \Holder-continuous on $(0,\infty)$, and $x^\alpha f(x)$ is \Holder-continuous at zero, for $\alpha < 1$.   Then the Cauchy transform  $\CC f(z)$ is the unique function $\phi(z)$ that satisfies all of the following properties:
\beginorderedlist
	\newitem analytic off $[0,\infty)$,
	\newitem has weaker than pole singularities,
	\newitem decays at infinity, 
	\newitem has continuous limiting values
	$$\phi_\pm(x) = \lim_{\epsilon \rightarrow 0^+} \phi(x \pm \I \epsilon) \qfor x\in (0,\infty),$$
and
	\newitem satisfies the jump condition
	$$\phi^+(x) - \phi^-(x) = f(x) \qfor x\in (0,\infty).$$
\endorderedlist
 Furthermore,
	$$\CC^+ f(x) + \CC^- f(x) = -\I \H f(x)\qfor x \in  (0,\infty).$$

\Remark
	This lemma is classical, with the precise conditions of this version combining \cite[Theorem 14.11a]{HenriciComplexAnalysis3} and \cite[Lemma 7.2.2]{FokasComplexVariables}.  Uniqueness follows from Liouville's theorem: the difference between any two solutions will be continuous on $(0,\infty)$ and have a singularity weaker than a pole at the origin, hence it is entire, and therefore zero due to decay at infinity.    Weaker variants of this lemma are available in terms of ${\rm L}^p$ spaces (cf.~\cite[pp.~100]{FokasPainleve}).   The derived formul\ae\ also carry over to the weaker setting.

Plemelj's lemma translates the task of calculating Cauchy transforms from an integral operator formulation to an analytic boundary value problem.  It also means that the Hilbert transform can be calculated directly from the Cauchy transform.    We will use Plemelj's lemma  to  determine the behaviour of the Cauchy transform under analytic maps.

	In \cite{SOHilbertTransform}, it was observed that a Cauchy transforms are in a sense invariant with respect to {\Mobius} transformations.  In our setting, this can be used to relate a Cauchy transform over $[0,\infty)$ to a Cauchy transform over $[-1,1]$:
\Proposition{mob}
	Let $M(x) = {x + 1 \over x - 1} L$ be a {\Mobius} transformation mapping  $[-1,1]$ to $[0,\infty)$.  Provided that $f$ satisfies the conditions of \lmref{plemelj}, then
	$$\CC f(z) = \CC_{(-1,1)}[f \circ M](M^{-1}(z)) - \CC_{(-1,1)}[f \circ M](1),$$
where
	$$\CC_{(-1,1)} f(z) = {1 \over 2 \pi \I} \int_{-1}^1 {f(x) \over x - z} \dx.$$
Therefore, 
	$$\H f(x) = \H_{(-1,1)}[f \circ M](M^{-1}(x)) -  \H_{(-1,1)}[f \circ M](1),$$
where
	$$\H_{(-1,1)} f(x) = {1 \over \pi} \dashint_{-1}^1  {f(t) \over t - x} \dt.$$

The validity of this proposition follows immediately from Plemelj's lemma.  The formul\ae\ derived below can be seen as generalizing this basic result to the case where the map $M$ is not conformal.

%
%

%

\Section{map} Monomial maps of Cauchy transforms.
	
Consider the map $x^p$ where $p$ is an integer.  We wish to use \propref{mob} to relate $\CC f$ to  $\CC[f(\diamond^p)]$, however, $x^p$ is not a {\Mobius} transformation, and has multiple inverses\footnote{$^\star$}{We use $\diamond$   in place of $\cdot$ to denote the independent variable.}.  The solution is to sum over {\it all} $p$ inverses.  The following proposition ensures that the resulting expression remains  analytic:

\Proposition{polyroots}
	Let $p$ be an integer and define the $p$ inverses of $z^p$ by  $\lambda_j(z) = z^{1/p} \E^{2 \pi \I j/p}$.  If  $g$ is an analytic function everywhere except on a closed set $\Sigma$,
then
	\mEq{gsum}{\sum_{j=0}^{p-1}  g(\lambda_j(z))}
is analytic in $z$ everywhere off $\Sigma^p := \set{z^p : z \in \Sigma}$.

\Proof

Suppose $z \neq 0$ is not in $\Sigma^p$, so that $\lambda_j(z) \notin \Sigma$ for all $k$.  Note that the value of \eq{gsum} does not depend on the choice of branch cut of  $z^{1/p}$: moving the branch cut  corresponds to merely a relabelling of $\lambda_j(z)$.    Thus we can choose a branch cut that avoids $z$, hence each term of $\eq{gsum}$ is analytic. If $0 \notin \Sigma$, analyticity at $z = 0$ follows because \eq{gsum} is analytic everywhere in a neighbourhood of $0$ and is continuous at zero.  

\mqed

We now obtain the following generalization:

\Lemma{Pol}
	Let $p$ be an integer and define  the $p$ inverses of $z^p$ by $\lambda_j(z) = z^{1/p} \E^{2 \I \pi  j/ p}$ .    If $f$ and $f(x^p)$ satisfy the conditions of \lmref{plemelj}, then the Cauchy transform satisfies
	\mEq{P}{\CC f(z) = \sum_{j=0}^{p-1} \CC[f(\diamond^p)](\lambda_j(z)).}
Therefore,
	$$\H f(x) = \H[f(\diamond^p)](x^{1/p}) + 2 \I \sum_{j=1}^{p-1} \CC[f(\diamond^p)](\lambda_j(x)).$$

\Remark
This result is a special case of a remark in \socite{SORHFramework; pp.~314}, which gives a formula for Cauchy transforms under general polynomial maps, but was stated without proof.  

\Proof

 We  show that
 	$\phi(z) =  \sum_{j=0}^{p-1} \CC[f(\diamond^p)](\lambda_j(z))$
  satisfies each of the conditions of Plemelj's lemma stated in \lmref{plemelj}:
%
\beginorderedlist
   \soitem $\phi$ is analytic off $[0,\infty)$ by the preceding proposition.    
   \soitem $\phi$ has  weaker than pole singularities since $\CC[f(\diamond^p)](\lambda_j(z))$ will have a weaker singularity at the origin than $\CC[f(\diamond^p)](z)$ itself.
   \soitem  The decaying property follows from the decay in each Cauchy transform.  
   \soitem Continuity of $\phi_\pm$ follows from continuity of each Cauchy transform and $\lambda_j(z)$.  
\soitem The jump condition is a consequence of applying Plemelj's lemma.  With the standard choice of branch cut and $x \in (0,\infty)$, $\lambda_0(x) = x^{1/p} \in (0,\infty)$, and we get, for $w(x) = f(x^p)$,  
	$$\br[\CC w(\lambda_0(x))]^+ - \br[\CC w(\lambda_0(x))]^- = \CC^+ w(x^{1/p}) - \CC^- w(x^{1/p}) = w(x^{1/p}) = f(x).$$
 For all other $j$, $\lambda_j(x)$ are bounded away from $[0,\infty)$, hence
	$\sum_{j=1}^{p-1} \CC w(\lambda_j(z))$
is analytic for $x \in [0,\infty)$, and has no jump.  We thus get
	$$\phi^+(x) - \phi^-(x) = \br[\CC w(\lambda_0(x))]^+ - \br[\CC w(\lambda_0(x))]^-  =  f(x).$$
\endorderedlist
The expression for the Hilbert transform follows immediately.

\mqed

\Figuretwofixed{Cubic}{Third}
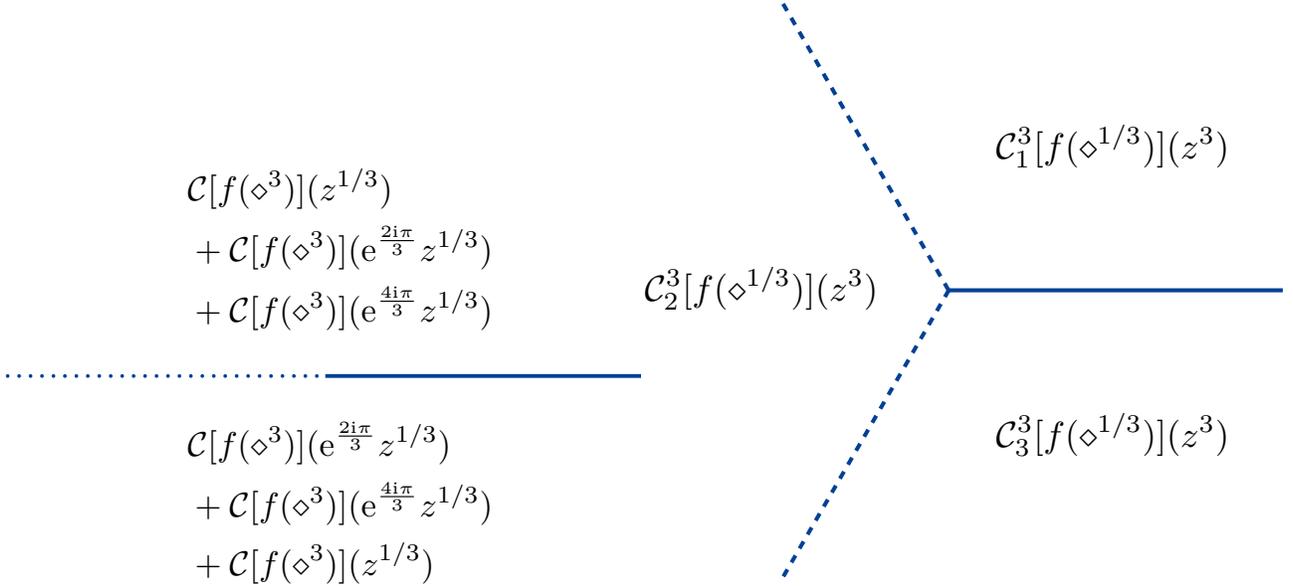
	$\CC f$  in terms of $\CC[f(\diamond^{3})]$ (left) and $\CC^3[f(\diamond^{1/3})]$ (right).  Solid line corresponds to the jump in $\CC f$, dotted line corresponds to the removable branch cut of $z^{1/3}$, and dashed line corresponds to the removable cut arising from changing sheets of $\CC^q$.

The particular formula for $p = 3$ is given on the left of  \figref{Cubic}.   Note that the representation with the standard branch cut for $z^{1/p}$ has a removable branch cut along $(-\infty,0]$.

\Section{inte}  Fractional monomial maps of Cauchy transforms and Riemann surfaces.

	  In this section, we investigate Cauchy transforms under the map $x^{1/q}$  where $q$ an integer.    
%
	  We begin with an example before stating the general result:

\Example{quad} 
	  Consider the special case of $q = 2$, so that we want to relate the Cauchy transform $\CC f$ to $\CC[f(\sqrt \diamond)]$.  If we attempt to use the  representation
	$$\CC f(z) \questionequals \CC[f(\sqrt \diamond)](z^2)$$
we encounter an issue: the right-hand side has a branch cut along negative $z$.  
	To overcome this, consider the solution to the following two scalar Riemann--Hilbert problems with $x \in [0,\infty)$,
	\meeq{
		\varphi^{(1)}_+(x) - \varphi^{(1)}_-(x) = f(\sqrt x) \qqand \varphi^{(1)}(\infty) = 0, \ccr
		\varphi^{(2)}_+(x) + \varphi^{(2)}_-(x) = f(\sqrt x) \qqand \varphi^{(2)}(\infty) = 0.
		}
By Plemelj's lemma, $\varphi^{(1)}$ is precisely $\CC[f(\sqrt \diamond)]$, while $\varphi^{(2)}$ can be constructed as $\P_{1/2}[f(\sqrt \diamond)]$, for $\P_r$  defined in \df{DVarphi}.
%
We claim that
	$$\Phi(z) = \half \socases{\varphi^{(1)}(z^2) + \varphi^{(2)}(z^2)  & 0 < \arg z < \pi \cr
						\varphi^{(1)}(z^2) - \varphi^{(2)}(z^2)  & -\pi \leq \arg z < 0}	 $$
satisfies the correct jump and decay properties.  We first note along $[0,\infty)$ we have
	$$\Phi^+(x) - \Phi^-(x) = {\varphi_+^{(1)}(x^2) + \varphi_+^{(2)}(x^2) - \varphi_-^{(1)}(x^2) + \varphi_-^{(2)}(x^2) \over 2}  = f(x).$$ 
While $\varphi^{(1)}(z^2)$ and $\varphi^{(2)}(z^2)$ have branch cuts along $[0,-\infty)$, $\Phi$ does not:
	$$\Phi^+(x) - \Phi^-(x) = \varphi_+^{(1)}(x) - \varphi_+^{(2)}(x) - \varphi_-^{(1)}(x) + \varphi_-^{(2)}(x) = 0.$$
Thus, since it has the same jump and decays at $\infty$,  $\Phi$ must be $\CC f(z)$ by Plemelj's lemma.

To extend the results, we first recall the  properties of the following integral transform:

\Definition{DVarphi}
	Define 
	$$\P_r f(z) = \E^{\pi \I r} (-z)^r \,\CC\!\br[f \over \diamond^r](z)$$ 
for $r \geq 0$.

	This operator satisfies the following analogue to Plemelj's lemma:

\Proposition{Varphi}
Assuming that $f x^{-r}$ satisfies the conditions of \lmref{plemelj},  then $\P_r f(z)$ satisfies the following:
\beginorderedlist
	\newitem analytic off $[0,\infty)$,
	\newitem decays at $\infty$, 
	\newitem has weaker than $\abs{z}^{r - 1}$ singularity at the origin,  
	\newitem has continuous limiting values $\P_r^\pm f(x)$ for $x \in (0,\infty)$, and
	\newitem satisfies the jumps
	\meeq{
		\P_r^+ f(x) - \E^{-2 \I \pi r} \P_r^- f(x) = f(x) \ccr
		\P_r^+ f(x) + \E^{-2 \I \pi r} \P_r^- f(x) = - \I x^r \H\!\br[f \over \diamond^r](x)
		}
\endorderedlist

\Remark
	The case $r = \half$ is a classical example, see e.g.~\socite{FokasComplexVariables; Example 7.3.3}, and the general result follows from the classical theory, see e.g.~\socite{SIE; Section 37}.   Alternatively, it is a direct consequence of Plemelj's lemma and the jump of $(-z)^r$.

We can use these operators to construct a generalization of the Cauchy transform living on a $q$-sheeted Riemann surface.    We represent the  Riemann surface by $q$ copies of the complex plane with branch cuts along $[0,\infty)$.  Then a function $\Phi(z) $ analytic on the Riemann surface is represented by  $q$  analytic functions $\Phi_\nu(z)$ satisfying the jumps
	$$\Phi_\nu^-(x) = \Phi_{\nu+1}^+(x)\qfor \nu = 1,\ldots, q - 1 \qqand \Phi_q^- = \Phi_1^+,$$
provided that $\Phi$ is analytic between each sheet.  

The generalization of the Cauchy transform 
takes the form of finding $\Phi(z)$ that satisfies the jump $\Phi^+(x) - \Phi^-(x) = f(x)$ on a {\it single sheet} of the Riemann surface, is analytic everywhere else on the Riemann surface, and decays at $\infty$ on {\it all sheets}.  We define this as follows, taking the jump to like in between the $q$th and first sheet:

\Definition{Cauchyrsf}
Define 	$\CC^q f(z)$ on the $q$-sheeted Riemann surface by
\meeq{
		q \CC_1^q f(z) =  \CC  f(z) + \P_r f(z)  + \P_{2r}  f(z)  + \cdots + \P_{1 - r}f(z) \ccr
		q \CC_2^q f(z) =  \CC  f(z) +  \E^{2 \pi \I \over q} \P_r f(z)  + \E^{4 \pi \I \over q}\P_{2r}  f(z)  + \cdots + \E^{2 {q - 1 \over q} \pi \I}\P_{1 - r}f(z) \ccr
		q \CC_3^q f(z) =  \CC  f(z) +  \E^{4 \pi \I  \over q} \P_r f(z)  + \E^{8 \pi \I \over q}\P_{2r}  f(z)  + \cdots + \E^{4 {q - 1 \over q} \pi \I}\P_{1 - r}f(z) \cr		
		&\vdots \ccr
		q \CC_q^q f(z) = \CC  f(z) +  \E^{2 {q - 1 \over q} \pi \I} \P_r f(z)  + \E^{4 {q - 1 \over q} \pi \I}\P_{2r}  f(z)  + \cdots + \E^{2 {(q - 1)^2 \over q} \pi \I}\P_{1 - r}f(z),
		}
for $r = 1/q$.

\Lemma{Cauchrsf}
 Assuming $f(x) x^{-{\nu \over q}}$  for $\nu = 0,\ldots,q-1$ satisfy the conditions of \lmref{plemelj}, $\CC_\nu^q f$ satisfy the following jump relationships for $x \in (0,\infty)$:
	$$\CC_\nu^{q-} f(x) = \CC_{\nu+1}^{q+} f(x)\qfor \nu = 1,\ldots, q - 1 \qqand \CC_1^{q+} f(x) -  \CC_1^{q-} f(x) = f.$$

\Remark
	Uniqueness is not needed in our setting, though see   \cite{CauchyRiemannSurface} for further discussion.  

\Proof
	For $\nu  = 1,\ldots, q - 1$ we have
	\meeq{
	q\CC_{\nu+1}^{q+} f - 	q\CC_{\nu}^{q-} f = \br[\CC^+ f - \CC^- f] + \E^{2 \nu \pi \I \over q} \br[\P_r^+ f - \E^{-2 \nu \pi \I \over q} \P_r^- f]+ \cdots \cr
	&\qquad+  \E^{2 \nu \pi \I {q - 1\over q}} \br[\P_{1-r}^+ f - \E^{-2 \nu \pi \I (1 -r)}\P_{1-r}^- f]\ccr
				=   \pr(1 + \E^{2 \nu \pi \I \over q} + \cdots +  \E^{2 \nu \pi \I {q - 1\over q}}) f 
 				=   {\pr(\E^{2 \nu \pi \I \over q})^q - 1 \over \E^{2 \nu \pi \I \over q} - 1} f = 0.
}
Similarly, 
	\meeq{
	q\CC_1^{q+} f- q\CC_q^{q-} f = \br[\CC^+ f - \CC^- f] + \br[\P_r^+ f - \E^{-2 \pi \I \over q} \P_r^- f]+ \cdots + \br[\P_{1-r}^+ f - \E^{-2 \pi \I (1 -r)}\P_{1-r}^- f] \ccr
						= q f.
}

\mqed

We can finally use this to reduce $\CC f$, by interpreting the map $z^q$ as traversing each sheet individually.   This can be written as follows, see also the right-hand side of \figref{Cubic}:

\Lemma{int}
Suppose  $q$ is an integer and  $f(x^{1/q}) x^{-{\nu \over q}}$  for $\nu = 0,\ldots,q-1$ satisfy the conditions of \lmref{plemelj}. Then
	$$\CC f(z)  = \CC_\nu^q [f(\diamond^{1/q})](z^q)  \qfor 2 \pi {\nu-1 \over q} < \arg z < 2 \pi {\nu \over q}.$$
  Therefore,
	\mEq{IntegerHil}{\H f(x) =  {1 \over q} \sum_{\nu=0}^{q-1} x^{\nu} \H\!\br[f(\diamond^{1/q}) \over \diamond^{\nu/q}](x^q). }

\Proof
	Analyticity away from $[0,\infty)$ and decay at infinity follow from the preceding lemma: we are simply moving from one sheet to another as $\arg z$ passes over $2 \pi r \nu$.  On $(0,\infty)$ we have the correct jump:
	$$\CC_1^{q+} [f(\diamond^{1/q})](x^q) - \CC_q^{q-} [f(\diamond^{1/q})](x^q) = f(x).$$
We only need to show that we have weaker than pole singularities.  For each $\P_{\nu/q} f$ in the definition of $\CC^q$, this follows directly from the third condition of \propref{Varphi}.  For $\CC[f(\diamond^{1/q})](z^q)$,  the fact that $x^{\alpha - {q-1 \over q}}  f(x)$ is \Holder-continuous at zero for  $\alpha < 1$ along with \socite{FokasComplexVariables; Lemma 7.2.2}  implies that we can bound the singularity at the origin:
	$$\abs{\CC[f(\diamond^{1/q})](z)} \sim {C \over z^{\alpha - 1 + {1 \over q}}}$$
The singularity  of $\CC[f(\diamond^{1/q})](z^q)$ is thus $z^{-q (\alpha - 1) - 1}$, which is  weaker than a pole since $\alpha  < 1$.  

\mqed

	\Remark For calculation of the related Hilbert transforms of the mapped function, we note that the left and right limits $\P_r^\pm f(x)$ can be calculated via \propref{mob} and numerical methods for the Hilbert transforms with Jacobi weights \cite{DellaHalfHilbert,GautschiJacobiHilbert,KimCauchyJacobi}.  For the special case of $r =1/2$, $\P_r f(z)$ can be calculated uniformly accurate throughout the complex plane  via modifying the procedure of \cite{SOHilbertTransform,SOEquilibriumMeasure} to derive a simple formula in terms of the Chebyshev  expansion of $f\fpr({x - 1 \over x + 1})$.   For general $r$, there does not appear to be a uniform numerical method, however, for $z$ away from the real axis  they can be reliably calculated by expanding $f\fpr({x - 1 \over x + 1})$ in Jacobi series, and using Olver's algorithm \cite{OlversAlgorithm} (see also \socite{DLMF;  \S3.5}) to calculate the Cauchy transforms of Jacobi polynomials, since they are the minimal solution of the Jacobi polynomial recurrence relationship \socite{GautschiOrthogonalPolynomials; \S2.3.1}.

\Section{rat} Rational monomial maps.

	We now turn our attention to the case of rational $r = p/q$. 
	We combine the two previous results: we  sum over  $p$  terms of the form  $\lambda_j(z) =  \E^{2 \pi \I  j / r} z^{1/r}$,  respecting the Riemann surface.   Two examples are given in \figref{TwoThirds}.  The following theorem states the general result:
	

\Figuretwofixed{TwoThirds}{ThreeHalfs}
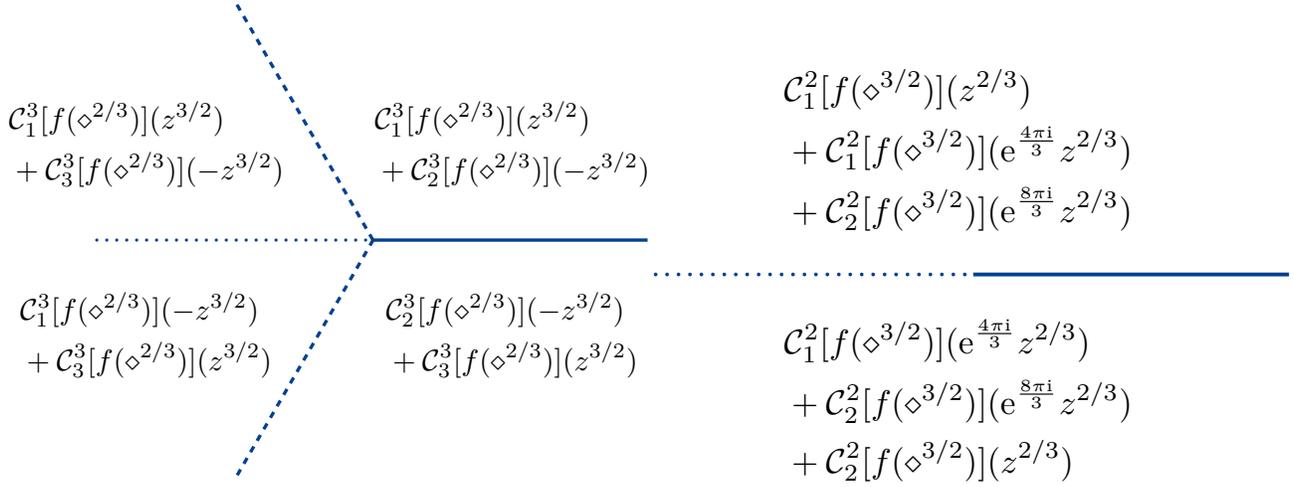
	$\CC f$ in terms of  $\CC^3[f(\diamond^{2/3})]$ (left) and $\CC^2[f(\diamond^{3/2})]$ (right).  Solid line corresponds to the jump in $\CC f$, dotted line corresponds to the removable branch cut of $z^r$, and dashed line corresponds to the removable cut arising from changing sheets of $\CC^q$.

\Theorem{rational}
Let $ r = p/q$ and $\lambda_j(z) = \E^{2 \pi \I j/r} z^{1/r}$.  Assuming that $f(x^r) x^{-\nu /q}$ for $\nu = 0,\ldots,q-1$ satisfy the conditions of \lmref{plemelj},
	$$\CC f(z) = \sum_{j = 0}^{p-1} \CC_{\beta(j,z)}^q[f(\diamond^r)](\lambda_j(z)),$$
where $\beta$ chooses the branch that $\lambda_j(z)$ lives on:
	 $$\beta(j,z) = \left\lceil {\arg z^{1/r} \over 2 \pi}  + j/r\right\rceil.$$
Therefore,
	$$\H f(x) = {1 \over q} \sum_{\nu=0}^{q-1} x^{\nu} \H\!\br[f(\diamond^r) \over \diamond^{\nu/q}](x^{1/r}) + 2 \I  \sum_{j=1}^{p-1} \CC_{\beta(j,z)}^q[f(\diamond^r)]( \lambda_j(z)).$$
	 
\Proof
	The proof  formula are derived by combining  \lmref{int} and \lmref{Pol}:
	\meeq{
		\CC f(z) = \sum_{j=0}^{p-1} \CC[f(\diamond^{p})](\tilde \lambda_j(z)) = \sum_{j=0}^{p-1} \CC_\nu^q[f(\diamond^{p/q})](\tilde \lambda_j(z)^q) \ccr
		= \sum_{j=0}^{p-1} \CC_{\beta(j,z)}^q[f(\diamond^r)]( \lambda_j(z))
		}
for $2 \pi {\nu-1 \over q} < \arg z < 2 \pi  {\nu \over q}$ and $\tilde \lambda_j(z) = \E^{2 \pi \I j/p} z^{1/p}$.     Similarly,
	\meeq{ \H f(z) = \H[f(\diamond^p)](x^{1/p}) + 2 \I \sum_{j=1}^{p-1} \CC[f(\diamond^p)](\tilde \lambda_j(x)) \ccr
				= {1 \over q} \sum_{\nu=0}^{q-1} x^{\nu} \H\!\br[f(\diamond^r) \over \diamond^{\nu/q}](x^{1/r}) + 2 \I  \sum_{j=1}^{p-1} \CC_{\beta(j,z)}^q[f(\diamond^r)]( \lambda_j(z)).
	}
The result can verified using \lmref{plemelj}, similar to preceding proofs.


%
\mqed


%
%
%

%
\Section{num}  Examples.

\Figuretwofixed{Changingr}{Optimalr}
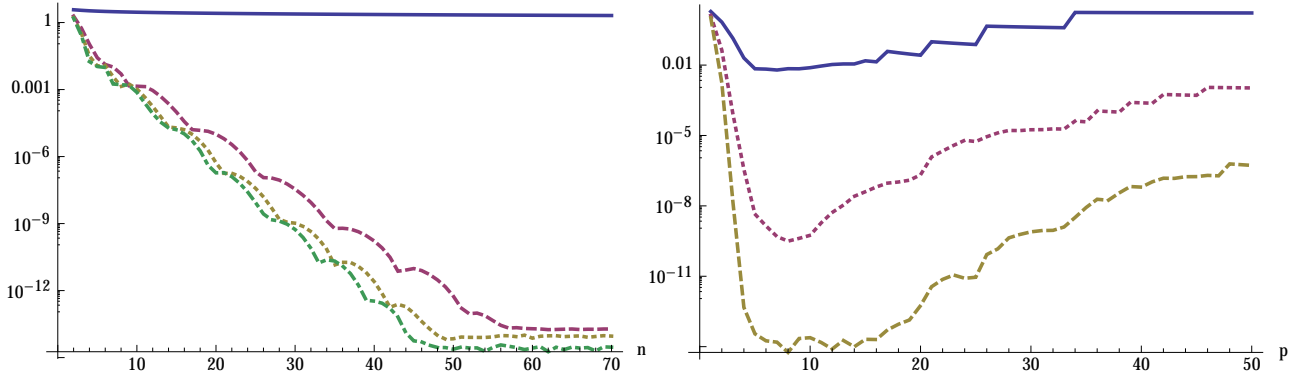
	The error of approximating $\CC f(1 + \I)$ for $f(x) = {x \over (1 + x)^{\pi -2}}$ by Cauchy transforms of $f(x^p)$, using \eq{pCauch}.  In the left graph, $p = 1$ (plain), 10 (dashed), 20 (dotted) and 100 (dot-dashed).  In the right graph, $n = \floor{100/p}$ (plain), $\floor{500/p}$ (dotted), and $\floor{1000/p}$ (dashed).

Our first example is one where $f$ has slow decay:
	$$f(x) = {x \over (1 + x)^{\pi -2}} \sim x^{-0.14159\ldots}.$$
We can choose $p$ large so that $f(x^p) \sim x^{-p 0.14159\ldots}$ has fast decay.  We truncate the interval into subintervals:
	  $$\vect[a_1,\dots,a_\ell] = \vect[0,10^{-10/p},1,10^{10/p},10^{20/p},\dots,10^{110/p}],$$ where the last endpoint $10^{110/p}$ is chosen  because $f(10^{110}) \approx 2.7\times10^{-16}$ is on the order of machine epsilon.    In each sub-interval $[a_l,a_{l+1}]$, we approximate
	 $$f\fpr({\pr({a_l + a_{l+1} \over 2} + {a_{l+1} - a_l \over 2} x)^p}) \approx \sum_{k = 0}^{n-1} f_k^l T_k(x),$$
where $T_k(x) = \cos k \arccos x$ are the Chebyshev polynomials and $f_k^l$ are calculated via the DCT.   We then  express the Cauchy transform in terms of Cauchy transforms of Chebyhev polynomials:
	\mEq{pCauch}{\CC f(z) \approx \sum_{j = 0}^{p-1} \sum_{k = 0}^{n-1} \sum_{l=1}^{\ell-1}  f_k^l \CC_{(-1,1)}T_k\fpr({2 \lambda_j(z) - a_l - a_{l+1}  \over a_{l+1} - a_l} ),}
where we used a finite-domain analogue of \propref{mob}  \socite{SORHFramework; Theorem 4.3} to reduce Cauchy transforms over $(a_l,a_{l+1})$ to $(-1,1)$.  
The Cauchy transform $\CC_{(-1,1)}T_k(z)$ are finally evaluating  in closed form via \socite{SOHilbertTransform; Theorem 6}.  In the left graph of \figref{Changingr}, we see that  exponential convergence to the true Cauchy transform is achieved as $n \rightarrow \infty$, and the rate of convergence improves slightly with increasing $p$.  However, the total computational cost is $\O(\ell n \log n + n p)$, and so in the right graph we compare choose $n$ depending on $p$ so that the computational cost is roughly fixed.  We see that, under this metric, $p \approx 10$ yields the highest accuracy.

\Figuretwofixed{ExponentialInteger}{SuperalgebraicInteger}
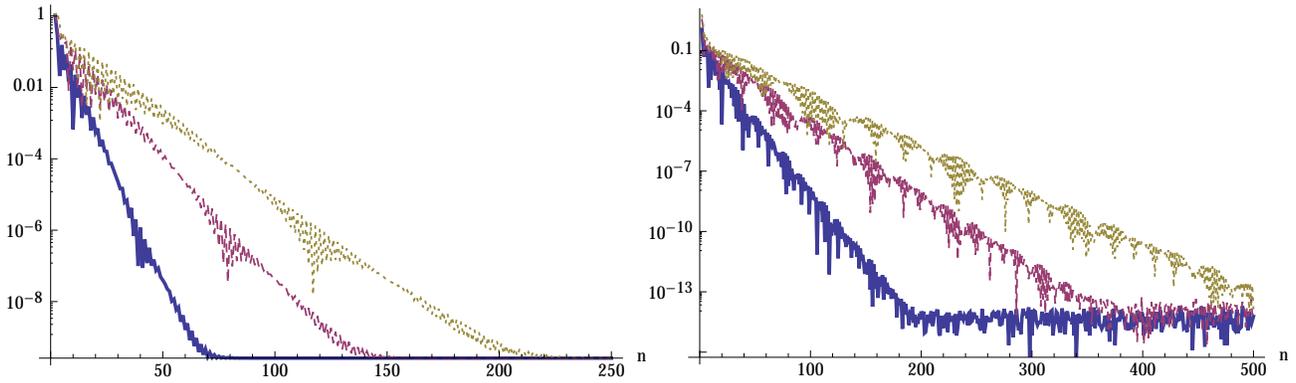
	The error of approximation $\H f(1.5)$ using for $f(x) = {x \over (1+x) (1+x^{1/5})}$ (left) and $f(x) = {1 + \E^{-x} \over 1 + x^{1/5}}$ (right), for $p = 5$ (plain), 10 (dashed), and 15 (dotted).  In the left figure, $f(x^p)$ has a converging series at $\infty$ and the approximation converges exponentially fast. In the right figure, $f(x^p)$ has only an asymptotic series and the approximation converges superalgebraically fast.

In our second example, we consider $\H f$ for the following two choices of $f(x)$ with decay like $x^{-1/5}$:
	$${x \over (1+x) (1+x^{1/5})} \qqand  {1 + \E^{-x} \over 1 + x^{1/5}}.$$
Under the map $x^p$ for $p = 5,10$, and 15,  $f(x^p)$ is smooth on $[0,\infty)$ and has a full asymptotic series at $\infty$.  We can thus employ the approach of \cite{SOHilbertTransform}.  Namely, we expand
	$$f\fpr({\pr({x + 1 \over x - 1})^p}) \approx \sum_{k = 0}^{n-1} f_k T_k(x)$$
where $f_k$ are calculated via the DCT.
 \propref{mob} is subsquently used to reduce the Hilbert transform:
	$$\H[f(\diamond^p)](x) \approx \sum_{k = 0}^{n-1} f_k \br[\H_{(-1,1)} T_k\fpr({x -1 \over x+1}) - \H_{(-1,1)} T_k\fpr(1) ].$$
We can thus caclulate $\H f$ by Combining this expression with the formula for the Hilbert transform of Chebyshev polynomials \socite{SOHilbertTransform; Theorem 6} and \lmref{Pol}.  
In  \figref{ExponentialInteger}, we see that the first example converges exponentially fast with $n$, since it has a converging asymptotic series at infinity, while the second example converges superalgebraically fast.  In both examples, increasing $p$ slows the convergence: the derivatives of the mapped function $f\fpr({\pr({x + 1 \over x - 1})^p})$ become increasingly large, slowing the convergence of the Chebyshev expansion.

For our third example, we  consider asymptotics of an oscillatory Hilbert transform,  with $f(x) =  \E^{(\I \omega  - 1) x^{3}}$.  
Near zero, we have
	$${ f(x^{1/3}) \over x^{\nu /3}} \sim {1 \over x^{\nu \over 3}} \sum_{k = 0}^\infty {(-)^k \over k!} x^{k} \E^{\I \omega x},$$
which implies that \socite{WangOscillatoryHilbert; Theorem 2.2}
	$$\H\br[f(\diamond^{1/3}) \over \diamond^{\nu/3}](x) \sim \I  {\E^{(\I \omega - 1) x} \over x^{\nu/3}} - {1 \over \pi} \sum_{\ell = 0}^\infty {\Gamma\fpr(\ell + 1 -{\nu \over 3})   \over \omega^{\ell + 1 - {\nu \over 3}}} \E^{{\pi \over 2} \pr(\ell + 1 - {\nu \over 3})\I} \sum_{j = 0}^\ell {(-1)^{j} \over j! x^{\ell - j + 1}}$$
as $\omega \rightarrow \infty$.  Using this expression for each term in \eq{IntegerHil} gives us an asymptotic expansion for $\H f(x)$.  In \figref{OscillatoryAsymptoticxHalf}, we compare the error as $\omega \rightarrow \infty$ for $x = \half$ and 2, taking $2, 4, 8$, and 16 terms in the asymptotic expansion.    We see that the approximation is indeed asymptotically accurate as $\omega \rightarrow \infty$, and also the accuracy improves with large $x$.

\Figuretwofixed{OscillatoryAsymptoticxHalf}{OscillatoryAsymptoticxTwo}
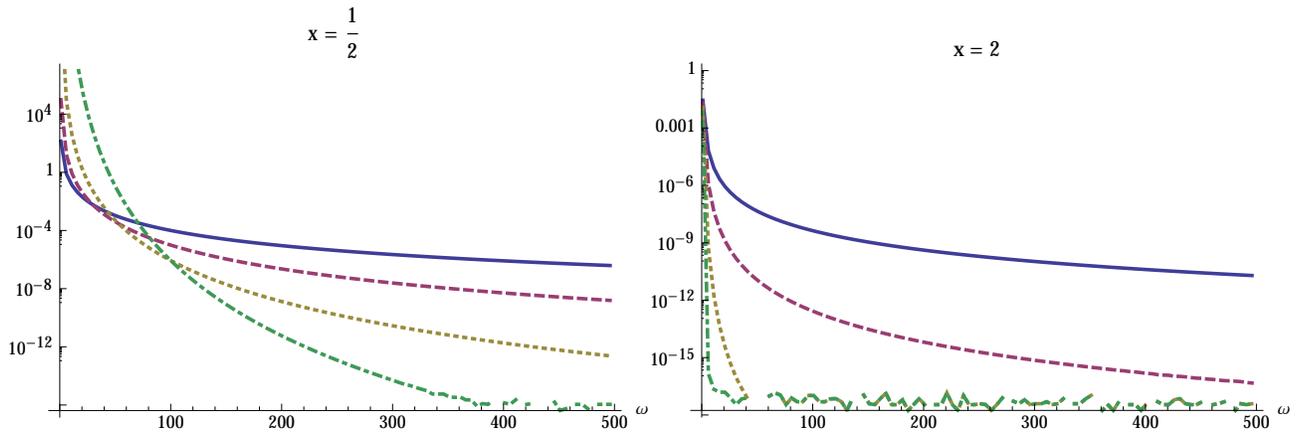
	The error in approximating $\H f(x)$ for $f(x) = \E^{(\I \omega + 1)x^3}$ by its $m$-term asymptotic expansion, for $m = 2$ (plain), 4 (dashed), 8 (dotted), and 16 (dot-dashed).  Observe that the approximations are better for large $\omega$ and large $x$.

\Acknowledgment
	I thank P. Antunes, whose question on numerical calculation of Hilbert transforms of functions with fractional decay motivated this research, T. Trogdon, for help with precise statement of Plemelj's lemma and helpful suggestions, and   A. Fish for interesting discussion of the results.  

\References

\Appendix{irr}    Irrational monomial maps.

If we approximate an irrational $r$ by a sequence of increasingly accurate rational numbers, \thref{rational} corresponds to an increasing number of summands and a generalized Cauchy transform living on a Riemann surface with an increasing number of sheets.  This suggests that the Cauchy transform of irrational $r$ can be expressed in terms of an {\it infinite sum} of evaluations of a  Cauchy transform living on an {\it infinite sheeted} Riemann surface.  We formally derive such a representation, and give numerical evidence supporting it.  

	The first step is to construct a generalized Cauchy transform living on an infinite number of sheets.  We do so formally by taking a limiting process of $\CC^q$.  On the first sheet we have
	\meeq{
		\CC_1^q f(z) = {1 \over q} \CC  f(z) + \P_r f(z)  + \P_{2r}  f(z)  + \cdots + \P_{1 - r}f(z) \ccr
		= \CC\br[ {1 + (z/\diamond)^{1/q} +  (z/\diamond)^{2/q} + \cdots + (z/\diamond)^{q - 1} \over q} f](z) \ccr
								= \CC\br[ {z - x \over q x [(z/x)^{1/q} - 1]} f](z)
		}
But as $q \rightarrow \infty$, we have $q \br[(z/x)^{1/q} - 1] \rightarrow \log {z \over x}$, which motivates the following definition for a generalized Cauchy transform $\CC^\infty$:
%
%
%
\Definition{Cinf}
	$$\CC_k^\infty f(z) = {1 \over 2 \pi \I} \int_0^\infty {f(x) \over x (\log x - \log(-z) - (2 k - 1) \pi \I)} \dx.$$
Numerically, this transform satisfies the following properties, justifying the notation:
\beginorderedlist
	\newitem decays at $\infty$ on every sheet,
	\newitem satisfies the jump 
	$$\CC_1^{\infty+} f(x) - \CC_0^{\infty-} f(x) = f(x),$$
	and
	\newitem  analytic between every other sheet 
	$$\CC_{\nu+1}^{\infty+} f(x) = \CC_\nu^{\infty-} f(x)\qfor k = 1,2,\ldots \qand k = \ldots,-2,-1.$$	
\endorderedlist

\Figuretwofixed{Irrational}{IrrationalPi}
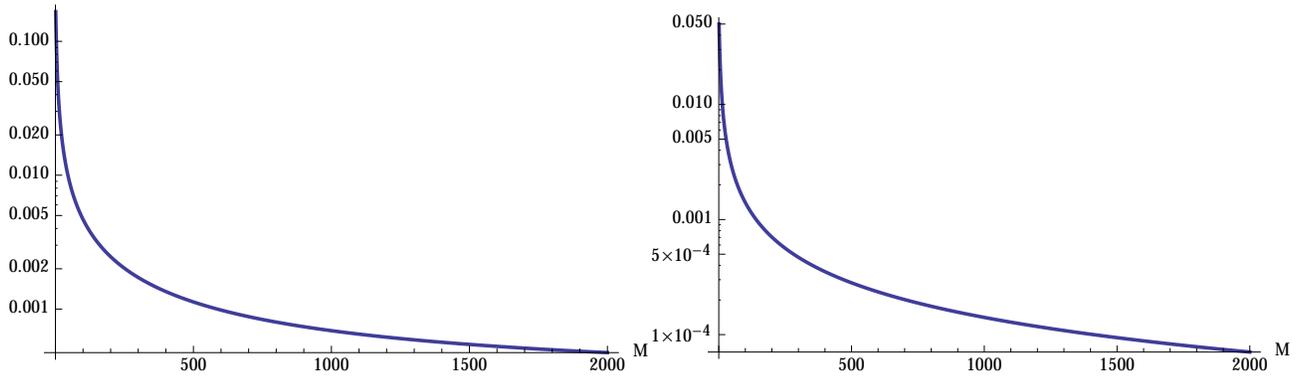
	For $r = \E$ (left) and $4-\pi$ (right), the error in approximating $\CC f(1+\I)$ by \eq{approxirr} for increasing $M$, where $f(x) = \E^{-x^{-1/r}} x^{-1/r}$.

We thus conjecture the following:
	$$\CC f(z) \questionequals \sum_{j = -\infty}^\infty \CC_{\beta(j,z)}^\infty[f(\diamond^r)](\lambda_j(z)) - \lim_{z \rightarrow \infty}  \sum_{j = -\infty}^\infty \CC_{\beta(j,z)}^\infty[f(\diamond^r)](\lambda_j(z)),$$
where $\beta$ and $\lambda_j$ are defined  in \thref{rational} and the second summation ensures that the representation decays at infinity.  We  can compare $\CC f(z)$ to the partial sum
	\mEq{approxirr}{\sum_{j = -M}^M \set{{\CC_{\beta(j,z)}^\infty[f(\diamond^r)](\lambda_j(z)) -   \CC_{\beta(j,z)}^\infty[f(\diamond^r)](\lambda_j(10^9 \I))}},}
where all transforms are performed numerically.  
In \figref{Irrational}, this approximation appear to converge as $M \rightarrow \infty$  for two choices of irrational $r$.  However, proving the accuracy of this representation appears delicate.  Furthermore, it is not clear how to efficiently numerically calculate the transform $\CC^\infty f$  and the slow convergence of the summation with $M$ means that it is of limited numerical use.

\ends